\input amstex
 \input amsppt.sty
 \documentstyle{amsppt}
 \loadmsbm
 \NoBlackBoxes
 %\pagewidth {4 in}
  %PREAMBLE
 \leftheadtext {(S. Liriano)}
 \rightheadtext{(Invariants of Parafree Groups)}
  \topmatter
 \title Algebraic Geometric Invariants
 of Parafree Groups
 \endtitle
 \author Sal Liriano
 \endauthor
 %\affil
 %\endaffil
 \address \flushpar
 % Dept. Of Science and Mathematics, B Building Room 831, FIT, Seventh
 %  Avenue  at 27th Street New York, NY 10001-5992
 \endaddress
 \curraddr
 \endcurraddr
 \email SAL21458\@yahoo.com
 \endemail
 \abstract

 Given a finitely generated ({\it  fg}) group $G$, the set
  $R(G)$ of homomorphisms from $G$ to $SL_2\Bbb C$
 inherits the
 structure  of  an  algebraic variety known as
 the {\it  representation variety \/}
 of $G$ in $SL_2 \Bbb C$. This algebraic variety is an invariant
 of {\it fg} presentations of $G$.
 Call a group $G$ parafree of rank $n$ if it shares the
 lower central sequence
 with a free group of rank $n$, and if it is residually nilpotent.
 The deviation of a {\it fg}
  parafree group
 is the difference between the minimum possible number of generators
 of $G$ and the rank of $G$.
 So  parafree groups
 of deviation zero are actually just free groups. Parafree groups that
 are not free share
 a host of properties with free groups. In  this paper
 algebraic geometric invariants involving the number of maximal
 irreducible components ({\it mirc}) of $R(G)$,
 and the dimension of $R(G)$
 for certain classes of
 parafree groups are computed.  It is
  shown that in an  infinite
 number of cases these invariants successfully discriminate
 between isomorphism types within
 the class of
 parafree groups of the same rank. This is quite
 surprising,
 since an $n$
 generated group $G$ is free of rank $n$ iff $Dim (R(G))=3n$. In
 fact, a direct consequence of Theorem 1.6 in this paper is that
 given an arbitrary positive integer $k$, and any integer $r \ge 2$
 there exist infinitely many non-isomorphic {\it fg}
 parafree groups of rank $r$ and deviation one
 with representation varieties of dimension $3r$, having more
 than $k$ {\it mirc} of dimension $3r$.
 This paper also
 introduces many novel and surprising dimension formulas for
 the representation varieties of certain one-relator groups.
 \endabstract
 \endtopmatter
 \document

\flushpar {\it General Structure of the Paper.} This paper begins
with an introduction where relevant ideas to what will follow are
developed. It then goes on to define what a parafree group is, and
how the notions that inspired G. Baumslag to give rise to such
groups arose in the context of investigations conducted by W.
Magnus, and a question of Hanna Neumann. The new results in this
paper are Theorem 1.0, Theorem 1.1, Theorem 1.2, Corollary 1.2,
Theorem 1.3, Theorem 1.5, and Theorem 1.6. In Section One several
results from the author's earlier work are introduced. In Section
Two the following theorems are proven: 1.0, 1.1, and 1.6. The
paper ends with a list notation.

 \head Introduction
 \endhead
 \medskip
The methods introduced here were developed with the explicit
agenda of exploring possible algebraic geometric parallels between
the representation varieties of a class of {\it fg} one relator
parafree groups and free groups. More precisely, the development
of this approach arose as a result of a question posed to the
author by G. Baumslag inquiring about the possible equality of the
Krull dimension of the space of $SL_2 \Bbb C$ representations of a
rank 2 free group, and that of the group $G=<x,y,z;x^2y^3=z^5>$, a
group that some readers may readily identify as a parafree group
of rank two, and deviation one.

\medskip
Let  $G$ be a group generated by the finite set
 $X = \{ x_1,x_2,\dots ,x_n \}$;
 then the set  $R(G)$ of homomorphisms
 from the
 group $G$ to $SL_2 \Bbb C$
 can  be  endowed  with  the  structure
 of a complex affine variety arising as the set of
 solutions in $ (SL_2 \Bbb C)^n$
 to the set of
 matrix
 equations obtained from the relations of a presentation of
 $G$ on the generators
 $X$. The affine variety $R(G)$  is an invariant of
 {\it fg} presentations of the group
 $G$. The invariant $R(G)$ exports into the study
 of {\it fg} groups the
 numerous invariants of {\it Algebraic Geometry\/}  and  {\it
 Commutative  Algebra\/} associated  with  algebraic  varieties.
 For a detailed account of the algebraic variety $R(G)$ the reader may
 consult \cite {LM}.

\medskip
 Given a  group $G$, invariants of  particular
 interest associated with $G$ are  the
 dimension, and  reducibility  status  of  $R(G)$. In the sequel,
 $Dim(R(G))$ shall denote the dimension of $R(G)$. Fortunately, if
 $ F_n$ is
 the free group
 of rank  $n$, then using the fact that $SL(2,\Bbb C)$ is an
 irreducible 3 dimensional variety and Lemma .6  it is
 easy to see that $R(F_n)= SL(2, \Bbb C)^n$, and is thus an
 irreducible
 \footnote {The product of irreducible varieties is irreducible.}
 variety of dimension $3n$.

 \medskip
 Denote by
 $\gamma_n G $ the $n$-th term of the  lower  central  series  of  the
 group $G$. A direct consequence of W. Magnus' 1935 paper
 \cite {MW}
 is what in the sequel will be referred to as Magnus' Theorem:

 \proclaim {Magnus' Theorem}
  A $k$-generated group $G$
 with
 $G/ \gamma_n G \cong F_k/ \gamma _n F_k $ for all $n\ge 1$ is free
 of rank $k$, where
 $F_k $ denotes the free group of rank $k$.
 \endproclaim

\flushpar
 Magnus' result led Hanna Neumann  to inquire whether it was  possible  for  two
 residually
 nilpotent groups $G$ and $G'$ to have
 $G/ \gamma_n G \cong G'/ \gamma _n G' $ for all $n$ without
 $G \cong G'.$
 Subsequently this question  prompted G. Baumslag   \cite {B1} to
 construct in the 1960's a class of
 groups he named {\it  parafree \/}.
 A group $G$ is termed
  parafree  if:
  \roster
 \item "{1)}" G is residually nilpotent.

 \item "{2)}" There exists a free group $F$ with the property that
 $G/ \gamma_n G \cong F/ \gamma _n F $ for all $n\ge  1$.
 \endroster

\flushpar
 Now denote by $\mu(G)$ the minimal number of generators of  a
 group  $G$.  Define the {\it  rank \/ } of $G$, denoted here
 by $rk(G)$, to be
 $rk(G) = \mu (G/ \gamma_2 G)$. Further, define the
 {\it  deviation  \/ }
  of the group $G$, here denoted by $\delta(G)$, to be
 $\delta(G) = \mu(G)-rk(G)$. A  parafree  group $G$ is
 of rank $r$ if the free group in (2) above  is also of rank $r$.
 Notice that by Magnus' Theorem  a parafree group of finite rank $n$ is
 free iff it has deviation zero.
 \medskip
 G. Baumslag  in  \cite { B1}
  introduced a result
 quite handy in building non-isomorphic parafree groups  of the
 same rank as a previously
 given parafree group:

 %Theorem 1, G. Baumslag from now on known as G. Baumslag's Theorem

  \proclaim { G. Baumslag's Theorem}
 Let $r$ and $n$ be positive integers and $H$ parafree of rank $r$, and
 let $(x)$ be the infinite cyclic group on $x$. Further, let
  $W\, \in \, H $. Suppose $W$ is the $k$ power of $W'$
 modulo $\gamma_2 H$, where $W'$ is itself not a power modulo
  $\gamma_2 H$; also, assume  that $k$ and $n$ are coprime
and that the generalized free product $ G = \{ H \ast (x) ; W =
x^n \}$
  is residually nilpotent. Then $G$ is parafree of rank $r$.
 \endproclaim

 \medskip
 Now consider the groups:
 $$G_{p_1\dots p_n} = \langle a_1,...,a_n;a_1^{p_1}a_2^{p_2}\dots
 a_{n-1}^{p_{n-1}} a_n ^{p_n}=1 \rangle, \tag 1-1 $$
 where $p_1,p_2,\dots,p_n$  are  positive  integers (all $\ge 2$ ).
 \medskip\flushpar
 In 1975, S. Meskin \cite {MS}
 showed that two groups $G$  and $G'$  as  in (1-1) are isomorphic iff
 there exists a permutation $\sigma \in \bold S(n)$ such that the sequence
 of exponents  $(p_1, p_2,...,p_n )$ of $G$ are sent by $\sigma$ to the
 corresponding sequence of  exponents  of  $G'$;  in  particular,  the  two
 sequences must be of the same length.
 \medskip
 \proclaim { Proposition 1}
 The groups $G_{p_1\dots p_n}$, with $n \ge 3$, are freely indecomposable
  parafree  groups of rank
 $(n-1)$ and deviation 1
 whenever the   $p_1,p_2,\dots,p_n$  are  positive
 integers with all the $p_{i's}\ge 2$ and
 having  no common  divisor.
 \endproclaim

 \demo {Proof}
  Proposition 1 is a direct result of G. Baumslag's Theorem  cited
  earlier, and
  Theorem 1
  in \cite {B6}, which guarantees that the groups are residually
  nilpotent. That each of the $G_{p_1\dots p_n}$ is of deviation one
  follows from the fact that $G_{p_1\dots p_n}$ maps to a corresponding
  free product $\Bbb Z_{p_1}* \cdots *
  \Bbb Z_{p_n}$, and an application of the Grushko-Neumann Theorem.
  Their free indecomposability can be deduced from Theorem 3 in
  \cite {SA}.
  \enddemo

 \medskip
 It is a well documented fact that one-relator cyclically pinched
 groups share
 a host of properties with free
 groups. For example, the fact that the groups
 in  (1-1) are cyclically pinched one-relator groups which meet the
 conditions of a combination theorem of \cite {BF} when $n\ge 3$ makes
 them
 hyperbolic, as is the
 case with {\it fg} free groups. In
 fact, when  $n\ge 3$ the groups
 in (1-1) are also linear by a result of Shalen \cite {SP}. These
 are all properties
 possessed by free groups. For additional properties shared by
 cyclically pinched one-relator groups and free groups, the
 reader may consult \cite {B4}, and  \cite {B5}.
 Incidentally,  as of this writing the existence on non-linear
 {\it fg} parafree groups is unknown; also unknown is  whether all
 finitely presented ones are hyperbolic.

  \medskip
  In \cite {L1} the author introduced a result
  permitting one
 to calculate the dimension and reducibility status of the
 representation variety of
 cyclically
 pinched one-relator groups obtained from a free group $F_n$ by
 adding a  new
 generator $y$ and a single relation $g=y^p$, for any non-trivial
 $g\in  F_n$. The use of the main result of \cite {L1} applied to
 the groups $G_{p_1,\dots,p_n}$, for $n \geq 3$, leads
 to Theorem 1.1, and many of the surprising
 results of this  paper.
 Before stating
 Theorem 1.1,  note that it takes into consideration
 the restrictions on the exponents of the relator in each of the groups
 $G_{p_1,\dots,p_n}$ in (1-1), and S. Meskin's result \cite {MS}.

 \medskip
 \proclaim {Theorem 1.1}
  Let $G_{pqt} =\langle  x_1,x_2,x_3; x^p_1 x^q _2 = x^t _3 \rangle $,
  where $p,q$ are negative integers $\le -2$, and $t\ge 2$ is an integer.
  Then,
  \roster
  \item "{i)}"  $Dim \,(R(G_{pqt})) = 6$,
  \item "{ii)}" $R(G_{pqt})$ is a reducible algebraic variety when at
  least one of the absolute values of $p, q, t$
   is strictly larger than two,
  \item "{iii)}" There exists an infinite set $S_2$ of groups
  associated with an infinite set $S$ of 3-tuples,  $S\subset \Bbb Z
  \times \Bbb Z \times \Bbb Z$,  having the property that if $(p,q,t)$,
  and $(p',q',t')$ are in $S$, then
    $R(G_{pqt}) \ncong R(G_{p'q't'})$
 if  $p \ne p'$,\, and  $q \ne q'$,\, and $t\ne t'.$
  \endroster
  \endproclaim

 \medskip \flushpar
 In fact, the next result gives a better picture of the general
 situation.

 \proclaim {Theorem 1.2} Let $p_1,p_2,\dots,p_n$ be positive  integers
 (all $\ge 2$ ). Then
 for $n \ge 3, \, Dim\, (R(G_{p_1, \dots , p_n})) = 3(n-1)$.
 \endproclaim

 \demo { Proof}
 Let $\Lambda_n$ be
 the $3(n-1)$ dimensional affine variety of Proposition 4 in Section
 Two. Clearly,
 $R(G_{p_1 \dots p_n}) \subset \Lambda_n$. So
 $Dim (R(G_{p_1\dots p_n})) \le  3(n-1)$. But, $G_{p_1 \dots p_n}$ has a
 presentation of deficiency $n-1$,  and
 hence $Dim\,(R(G_{p_1 \dots p_n })) \ge 3(n-1)$, by \cite {LM}.
 So $Dim\,(R(G_{p_1 \dots p_n}))=3(n-1)$.
 \enddemo

 \proclaim {Corollary 1.2}
 For any integer $r \ge  2$  there  exist  an infinite number  of  freely
 indecomposable pairwise non-isomorphic parafree
 groups of rank $r$ and deviation one  with $Dim \, (R(G)) =3r.$
 \endproclaim

 \demo {Proof} The proof follows from Theorem 1.2
 and Proposition 1.
  \enddemo

 Notice that Corollary 1.2 gives no information about the isomorphism
 type of the corresponding algebraic varieties. So it is conceivable
 that two non-isomorphic groups could have isomorphic representation
 varieties. In fact, in contrast to the case $n=3$, it is not known
 when, for fixed $n\ge 4$, and $p_1, \dots , p_n $ as in (1-1), the
 corresponding algebraic varieties $R(G_{p_1, \dots , p_n})$
 are reducible.

 \medskip\
 The next theorem, a direct consequence of Theorem 1.1,
 provides a rather strong illustration of the sensitivity of the
 invariant $R(G)$.
 \medskip
 \proclaim {Theorem 1.3}
 There exists an infinite set $S_2$ of freely indecomposable
  parafree  groups of rank two and deviation one having the
  property that for different
 $G_{pqt}$ and $G_{p'q't'}$ in $S_2$ their
 corresponding representation
 varieties are reducible, six dimensional, and non-isomorphic.
 \endproclaim

 \demo {Proof}
 Let $S_2$ be the infinite set of groups guaranteed by Theorem 1.1
 part iii. By Proposition 1 each of the groups in $S_2$ is
 parafree of deviation 1, rank 2, freely indecomposable,
 and for different $G_{pqt}$ and
 $G_{p'q't'}$ in $S_2$
 it is the case that $R(G_{pqt})\ncong R(G_{p'q't'}).$
 \enddemo

 \medskip\flushpar
 Theorem 1.3 displays clear
 distinctions in the reducibility
 status of the representation variety of  a rank 2 free group,
 and the
 representation varieties of an infinite number
 of non-isomorphic parafree groups all of rank 2 and deviation
 one. In fact all the groups in $S_2$ have non-isomorphic
 varieties, quite a difference from the free group of rank $n$,
 where up to isomorphism there is only one representation variety.
 The next theorem uses an invariant of $R(G)$ to characterize as
 free, or not free,
 an $n$-generated group.

  \proclaim {Theorem 1.4} An $n$ generated group $G$ is free of
  rank $n$ iff $Dim(R(G))=3n$.
  \endproclaim

  \demo {Proof} Let $G$ be $n$-generated with
  $Dim(R(G))=3n$; and suppose $G$ is not free. Then $G$ is a
  quotient group of $F_n$. So there exists a non-trivial word $w\ne 1$
  in $F_n$ such that $w=1$ in $G$. But by
  Theorem 1 in \cite {L2} one gets\footnote {$N(w)$ is the
  normal closure of $w$ in $F_n$.}  that $Dim(R(F/N(w)))< 3n$. But,
  $R(G)$ injects into $R(F/N(w))$, and consequently $Dim(R(G))\le 3n-1$.
  This is a contradiction. Conversely, suppose that $G$ is free and generated by $n$
  elements. Then by a result of Nielsen proven also in \cite {L2},
  these $n$ elements generate
  $G$ freely. Thus $Dim (R(G))=3n$. The proof is complete.
  \enddemo

 \medskip
 Let $V$ be an algebraic
 variety and let $N_c(V)$ be the integer
 valued function that counts its number of maximal irreducible
 components ({\it mirc})
 of dimension exactly $c$. Clearly $N_{c}((R(F_n))= 1$, when $c=3n$, and is
 zero otherwise. So that the reader gains a glimpse of the sensitivity of this
 invariant of algebraic varieties, and thus of {\it fg}
 groups, consider the case where $\frak T$ is a torus knot group. In
 \cite {L3} the author showed that $N_4(R(\frak T))=\bold g$, where $\bold g$ is
 the genus of the
 torus knot corresponding to the group $\frak T$.  The next theorem  is a
 consequence of ideas employed in the proof of Theorem 1.1
 involving the {\it mirc} counting function $N_c(V)$.

 \proclaim {Theorem 1.5}
  For any integer $r\ge 2$ there exists an infinite set
 $S_r$ of  parafree groups of rank $r$ and
 deviation 1 having the property  that for each group $G_i \in
 S_r$, $Dim(R(G_i))=3r$, and such that
 given another $G_j$ in $S_r$, then
 $R(G_i) \cong R(G_j)$  iff $G_i
 \cong G_j$.
 \endproclaim

 \demo {Proof}
 If $r=2$, then let $S_r=S_2$, as in Theorem 1.3. If $r \ge
 3$ let $S_r =\{F_{r-2}\ast G_i \vert G_i \in S_2 \}$,
  where
 $F_{r-2}$ is the free group of rank $r-2$. Note that any $G_k
 \in S_r$
 maps onto  $F_{r-2}\ast \Bbb Z_p  \ast \Bbb Z_q  \ast \Bbb Z_t $,
 where
 $p,q,t$ are all greater than 1. It follows from the Grushko-Neumann
Theorem that $G_k$ can't be generated by fewer than $r+1$
generators. That $G_k$ is parafree of
 rank $r$ follows from \cite {B2}, which guarantees that the rank of a
 free product of parafree groups is parafree of rank the sum of the
 factors. The dimension condition for $R(G_k)$ follows from
 Lemma .6, Theorem 1.1 and the fact that
 $Dim(R(F_{n-2}))=3(n-2)$. The
 isomorphism condition for different groups in $S_r$
 follows since any two groups  $G_e$, $G_f$ in $S_2$ have
 the property that $N_6(R(G_e)) \ne N_6(R(G_f))$, unless $G_e$ and
  $G_f$ are isomorphic, (see proof of Theorem 1.1 part iii); now by
  virtue of Lemma .2
  the same condition can be extended to the groups in $S_r$ by
  employing the {\it mirc } counting function $N_{3r}(V)$.
  The proof is complete.
 \enddemo

 The next theorem is perhaps one of the most surprising results
 so far on the representation varieties of groups in the class of
 {\it fg} parafree groups of rank $n\ge 2$, and deviation one.

 \proclaim {Thorem 1.6}
 Given an integer $r\ge 2$, and an arbitrary integer $k\ge 1$,
 there exists a parafree group $G$ of rank $r$ and deviation one
 with $Dim(R(G))=3r$, and  $N_{3r}(R(G)) \ge k$. If $r=2$,
 the parafree group can be taken to be freely indecomposable.
 \endproclaim

 \medskip Compare the above result with the situation in
 the case of parafree groups of rank $n$ and deviation
 zero. In such a case $N_{3n}(R(G))=1$, and up to isomorphism there
 is only one such parafree group, namely $F_n$.

 \medskip
 \head  Section One
 \endhead
 In this section Theorem .1 of \cite {L1} will be introduced
 along with some preliminary material to be employed in the
 eventual proof of many of the new results of this paper.

 \medskip
 Let $\,\,W \ne 1\,\,$
 be a freely reduced word in $F_n$ involving all the generators
 $\{ x_1, \dots, x_n \}$ of $F_n$. To the free group $F_n$ add a new
 generator $y$, and now consider the relation $W=y^k$  of the
 one-relator group
 $$G=\langle  x_1,\dots,x_n, y;W=y^k \rangle,\tag 2-1.1 $$
 where $k\ge 2$ is a positive integer.
 \medskip
 \remark {Observation 1 }
 The  relation  $W=y^k$ gives rise to an equation in
  $SL_2 \Bbb C$.
 Solutions  to  this  equation  are  $(n+1)$-tuples    of
 $SL_2  \Bbb C$ matrices
 $(m_1,\dots,m_{n+1})$  such   that   the   relation  $W=y^k$ is
 satisfied   when  the   $(n+1)$-tuple is  evaluated   in  $W=y^k$ under
 the obvious
 assignment $x_i \rightarrow m_i$, for all $x_i \in \{ x_1,  \dots,  x_n
 \}$, and $y\rightarrow m_{n+1}$. Further, observe that
  the relation can also be used to give rise to  an
 equation $W=-y^k$ in $SL_2 \Bbb C$.
 \endremark
 \medskip
 \remark {Notation }
 Denote the  $SL_2 \Bbb C$ solutions to the
 equations  $W=y^k$, and $W=-y^k$  of the previous observation  by
 $\Sigma_n( W, k)$, and  $\Sigma_n(-W, k)$ respectively.
 By $I$ shall  be  meant  the
 $2\times 2$ identity matrix.
 \endremark
 \medskip
 \flushpar Notice that
 $\Sigma_n(W,k)$ and  $\Sigma_n(-W,k)$  are
 affine algebraic  varieties, also  that the representation variety $R(G)$
 of the one-relator group in (2-1.1) is $\Sigma_n(W,k)$.
 Now  introduce the following integer valued function:
 $$ f(x)=\left\{ \aligned  1 \, \text{if} \,  x  &\ne  2\\  0
 \, \text {if} \, x&=2
 \endaligned \right. \tag 2-1.2 $$
 and let $P_+$ and $P_-$ be the  two  algebraic  sub-varieties  of  $R(F_n)$
 given by:
 $$P_+ =\{\rho \mid \rho \in  R(F_n),\text  {and}\,  \rho
 (W)=I \},  \qquad
 P_-=\{\rho\mid \rho \in R(F_n),\text {and}\, \rho(W)=-I\}. \tag 2-1.3$$
 \medskip
 \proclaim {Theorem .1}
 \roster
 \item "{a)}" $ Dim\,(\Sigma_n(W,k))=Max\{ Dim\,(P_+)+2f(k),
 Dim\,(P_-)+ 2, 3n \}\le 3n+1.$
 \item "{b)}" $Dim\,(\Sigma_n(-W,k))=Max\{Dim\,(P_-)+2f(k),
 Dim\,(P_+)+2, 3n\}\le 3n+1.$
 \endroster
 \endproclaim
\flushpar For a proof consult \cite {L1}.

\medskip
The next  corollary, a consequence of the fact that an irreducible
algebraic variety can have  no proper sub-variety  of its
dimension, will be indispensable.

 \proclaim {Corollary .1}
 \roster
 \item "{i)}"
 If in Theorem .1a  $Max\{Dim\,(P_+)+2f(k), Dim\,(P_-))+ 2 \}\ge 3n$,
 then $\Sigma_n(W,k)$ is a reducible variety.
 \item "{ii)}"
 If in Theorem .1b
 $Max\{Dim\,(P_-)+2f(k), Dim\,(P_+)+2\}\ge 3n$, then  $\Sigma_n(-W,k)$ is
 a reducible variety.
 \endroster
 \endproclaim
\flushpar For a proof consult \cite {L1}.

 \medskip
 Now a basic lemma concerning
 $SL_2 \Bbb C$ solutions of
 equations of the type $x^p=\pm I$ will be given, but first some notation.
 \remark {Notation }
 Given a positive integer $p$ and $M \in SL_2 \Bbb C$,
 denote by
 $\Omega(p,M)$ the following set
  $\{A  \mid  A\in  SL_2 \Bbb C,\,\,  A^p=M\}$.
 Notice  that  $\Omega(p,M)$  is  an  algebraic variety.
 \endremark
 \medskip
 \proclaim {Lemma .1}
 \roster
 \item  "{i)}"  If  $p=2$, then  $Dim\,\Omega(p,I)=0$.
 \item "{ii)}" If $p>2$, then  $Dim\,\Omega(p,I)=2$.
 \item "{iii)}" For $p\ge2$, $Dim\, \Omega(p,-I)=2$.
 \endroster
 \endproclaim
 \flushpar The proof of Lemma .1
 can be found \footnote {In $SL_2 \Bbb C$ a matrix of finite order
 is conjugate to a diagonal matrix. } directly, or by consulting
 \cite {L1}.
%Lemma a became known as Lemma .2
 \proclaim {Lemma .2} Let $V$ be an irreducible variety with
 $Dim(V)=m$, and let $W$ be a reducible variety with $N_q(W)=c$. Then
 $N_{m+q}(V\times W)=c.$
 \endproclaim
 \demo {Proof} The dimension of a product of two
 varieties  is the sum of the dimension of each factor, and the
 number of {\it mirc} in the product is the product of
 the {\it mirc} in each factor. Now,
 Since there is only one {\it mirc} of dimension $m$ in $V$
 given the stated
 irreducibility of the algebraic variety $V$, and since
 there are a number $c$ of $q$-dimensional {\it mirc} in $W$,
 it follows that the number of $(m+q)$-dimensional {\it mirc}
 of $V\times W$ is also $c$.
 \enddemo

 \medskip The first application of Lemma .1 will be the proof of
 the next result.

 \proclaim {Lemma .3}
 By $\Sigma_1( -x^p_1, -q)$ denote the solutions over
  $SL_2 \Bbb C$ to the
 matrix equation  $x^p_1 x^q_2 =-I$,  where $p$ and
 $q$ carry the same signs and
 assumptions as in  Theorem 1.1.
 Then:
 \roster
 \item "{i)}" $Dim\,(\Sigma_1( -x^p_1, -q))=4$, if $|p|>2$ or $|q|>2$,
 \item "{ii)}" $Dim\,(\Sigma_1( -x^p_1, -q))=3$, if $|p|=2$ and $|q|=2$.
 \endroster
 \endproclaim
 \medskip
 \demo {Proof i)}
 It has been assumed that $p$ and $q$ in
 $x^p_1x^q_2 = -I$  are  negative.  Consider the
 equation
 $x^p_1=-x^{-q}_2$. Let
 $P_-=\{\rho|\rho \in R(F_1), -\rho(x^p_1)= I\}$ and
 $ P_+=\{\rho|\rho \in R(F_1), -\rho(x^p_1)=-I \}$.
 Notice that if $|p|>2$, then by Lemma .1
 $Dim\, (P_-)=2$, and by
 Lemma .1 above
  $Dim\,(P_+)=2$. By Theorem .1b
 $Dim\,(\Sigma_1( -x^p_1, -q))=
  Max\, \{ Dim\,(P_-)+2f(-q), Dim\,(P_+)+2, 3\} \le 3+1 $, where
 $f$ is the function in the statement of Theorem .1.
 Thus $Dim\,(\Sigma_1( -x^p_1, -q))$ $=4$.
 If on the other hand $|p|=2$, then by Lemma .1
 $Dim\,(P_+)=0$, and $Dim\, (P_-)=2$.
 Thus
 $Dim(\Sigma_1(-x^p_1, -q))=4$, since $-q>2$.
 This completes the proof of (i).
 \enddemo
 \medskip
 \demo {Proof ii)}
 Suppose $p$ has absolute value equal to two.
 Notice that
 $$P_-=\{\rho|\rho \in R(F_1), -\rho(x^p_1)= I\},
 P_+=\{\rho | \rho \in R(F_1), -\rho (x^p_1)=-I\}.$$
 By Lemma .1 one gets that
 $Dim\,(P_-)=2$, and that
 $Dim(P_+)=0$. Thus by Theorem .1b  it follows that
 $Dim\,(\Sigma_1( -x^p_1, -q))=
 Max \,\{ Dim\,(P_-)+2f(-q), Dim\,(P_+ )+2, 3\} \le 3+1.$
 So $Dim\,(\Sigma_1( -x^p_1, -q))=3$, since $-q=2$. This completes
 the proof of (ii).
 \enddemo

 \proclaim {Lemma .4}
 \flushpar Let $G$ be a {\it fg} group,
  and let $N$ be a normal subgroup of $G$. Let $N_{c}(V)$ be
 the {\it mirc} counting function in dimension $c$, and
 suppose that $Dim (R(G))=c$. Then:
 \roster
 \item "{i)}" $Dim(R(G/N))\leq Dim(R(G)), $
 \item "{ii)}"
 $N_c (R(G/N))\leq N_c (R(G)).$
 \endroster
 \endproclaim

 \demo {Proof}
 \flushpar i). $N$ is the normal closure of at least one non-trivial
 element from $G$. So $R(G/N)$ is the zero locus in $R(G)$ of at
 least one polynomial belonging to the coordinate algebra of $R(G)$. It
 follows that $Dim(R(G/N))\leq Dim(R(G)) $.
 \flushpar ii). By i), $Dim(R(G/N))\leq Dim(R(G)).$ Assume that
 $N_c(R(G/N))> N_c(R(G))$. Then there is an irreducible component
 containing a component of its own dimension. This contradicts the
 properties of finite dimensional noetherian topological spaces.
\enddemo
\medskip
 \proclaim {Lemma .5}
 Denote by $\Sigma_2(-x_1^p x_2^q, t)$
 the solutions  over  $SL_2 \Bbb C$  to  the matrix
 equation  $x^p_1 x^q _2 =- x^t _3 $,
 where $p, q, t$ carry the same signs and assumptions  as  in
 Theorem 1.1. Then  $Dim\,( \Sigma_2(-x_1^p x_2^q, t))=6$.
 \endproclaim
 \medskip
 \demo {Proof}
 Theorem .1b will be employed. Notice that
 $$ P_-=\{\rho|\rho \in R(F_1),\,  -\rho(x^p_1x^q_2)= I\},
 P_+=\{\rho |\rho \in R(F_1), \, -\rho(x^p_1x^q_2)=-I \}.$$
 By Lemma .3,
 $Dim(P_-)= 4\, \text {or}\, 3$, depending on
 $|p|$ and $|q|$. Now consider
 $P_+$. By Theorem .1a-b and Lemma .1 follows that $Dim\,(P_+)=
 4$.  Now by Theorem .1b
 $$ Dim(\Sigma_2(-x_1^p x_2^q, t)) =
 Max\,\{ Dim(P_-)+2f(t), Dim(P_+)+2, 6\} \le 6+1 $$
 and
 thus $Dim(\Sigma_2(-x_1^p x_2^q,t))=6$.
 This completes the proof of Lemma .5.
 \enddemo

 \proclaim {Lemma .6} Let $G_1$ and $G_2$ be {\it fg}
 groups,
 and $V_1$,\,\, $V_2$ be algebraic varieties;
 then $R(G_1 \ast G_2)=R(G_1)\times R(G_2)$, and
 $Dim (V_1 \times V_2)=Dim(V_1)+Dim(V_2).$
 \endproclaim

 \flushpar The proof of the above can be easily deduced
 from elementary facts found in \cite {MD}, and is left to the reader.

\medskip
 Next is introduced a result allowing computation of the
 dimension of certain algebraic varieties using fibers of regular maps.
 \proclaim {Proposition 2 }
 Let  $\phi:V\rightarrow W$ be  a  regular map between  two
 algebraic varieties,
 where  $W$  is
 irreducible and $Dim\,(W)=n > 0$.
 Let $V_1$ and $W_1$  be  two
 proper
 closed subvarieties of $V$ and $W$, respectively,
 such that the restricted
 map $\phi:V^\circ \rightarrow\, W^\circ$, where
 $V^\circ=V-V_1$ and $W^\circ=W-W_1$, is such that:
 \roster
 \item "{1)}" $ \phi:V^\circ \rightarrow W^\circ $ is onto.
 \item "{2)}" $\phi$ has zero dimensional  fiber  above  each  point  of
 $W^\circ$.
 \item "{3)}" $\phi ^{-1}(W^\circ)=V^\circ$.
 \endroster
 Then
 $Dim\,(Cl(W^\circ))=Dim\,(Cl(V^\circ))= n$, where
  $Cl(W^\circ)$ denotes the Zariski closure of $W^\circ$.
 \endproclaim
 \flushpar
 For  a  proof  the  reader  may  consult  \cite
 {L1}, or use elementary facts from \cite {MD}.

 \head Section Two
 \endhead
\medskip The next result, instrumental in the proof of Theorem 1.1,
and Theorem 1.6,  guarantees that an infinite sequence of pairs
$(G_i, \frak N_i)$ consisting of a group $G_i$ and a corresponding
normal subgroup $\frak N_i$
 will have an infinite subsequence  provided the {\it mirc} counting
 function $N_c(R(G_i/\frak N_i))$ achieves an infinite number of
 values as $i$ varies. In
such an instance, the corresponding subsequence can be used to
obtain a well ordered infinite set of groups using the {\it mirc}
counting function $N_c(V)$ on their corresponding representation
varieties.

\proclaim {Theorem 1.0} \flushpar
 Let $S_{i=1}^{\infty}(G_i,\frak N_i)$ be an infinite
 sequence of pairs $(G_1,\frak N_1),
(G_2,\frak N_2), \cdots$ consisting of {\it fg} groups $G_i$ and a
corresponding normal subgroup $\frak N_i$ of $G_i$. Let $N_c(V)$
be the {\it mirc} counting function in dimension $c$. Suppose that
$Dim(R(G_i))=Dim(R(G_i/\frak N_i))=c$, and that the set
$\tilde{S}=\{ N_c(R(G_j/\frak N_j))\,\vert \,(G_j,\frak N_j) \in
S_{i=1} ^{\infty}(G_i,\frak N_i)\}$ contains an infinite set of
integer points. Then: \roster
\item "{i)}" $S_{i=1}^{\infty}(G_i,\frak N_i)$ has an infinite
subsequence $S^2$ with the property that given two different pairs
$(G_j,\frak N_j),\, (G_k,\frak N_k)$ in $S^2$ then
$N_c(R(G_j))\neq N_c(R(G_k)),$
\item "{ii)}" For different pairs
$(G_j,\frak N_j),\, (G_k,\frak N_k)$ in $S^2$ then $R(G_j)\ncong
R(G_k),$
\item "{iii)}" The set of groups
$S_2=\{G \vert  \text {$G$ occurs in some term of $S^2$} \}$ can
be well ordered using the function $N_c(V)$.

\endroster
\endproclaim
\demo {Proof of part ii)} Since $N_c(V)$ is an invariant of a
finite dimensional algebraic variety $V$, the proof of part ii
follows from the proof of part i.
\medskip
\demo {Proof of part i)}

By assumption the set $\tilde{S}$ is infinite. By Lemma .4 ii, it
follows that since each element of $\tilde{S}$ is bounded above by
an element of the set of integer points
 $\tilde{B}=\{N_c(R(G_j))\,\vert
\,(G_j,\frak N_j) \in S_{i=1} ^{\infty}(G_i,\frak N_i)\}$ that the
set $\tilde {B}$ is also infinite, since no finite set of positive
integers can be an upper-bound for an infinite set of positive
integers. Thus the set $\tilde{B}^2=\{
(N_c(R(G_j)),N_c(R(G_j/\frak N_j)))\,\vert \,(G_j,\frak N_j) \in
S_{i=1} ^{\infty}(G_i,\frak N_i)\}$ is an infinite subset of the
lattice $\Bbb Z ^2$ with strictly positive entries. Use the set
$\tilde{B}^2$ to obtain the subsequence $S^2$ of $
S_{i=1}^{\infty}(G_i,\frak N_i)$ meeting the stipulated
requirements.
\medskip
\demo {Proof of part iii)}  Given $G_j$ and $G_k$ in $S_2$, by
virtue of i, it is the case that either $N_c(R(G_j))<
N_c(R(G_k))$, or $N_c(R(G_k))< N_c(R(G_j))$, or $N_c(R(G_k))=
N_c(R(G_j))$\, in which case $G_j\cong G_k$.
\enddemo
\enddemo
\enddemo
\medskip
 The next theorem is  a cornerstone in the
 demonstration of many subsequent results.
 \proclaim {Theorem 1.1}
 Let $G_{pqt} =\langle  x_1,x_2,x_3; x^p_1 x^q _2 = x^t _3 \rangle $,
 where $p,q$ are negative integers $\le -2$, and  $t\ge 2$ is an
 integer.
 Then,
 \roster
 \item "{i)}"  $Dim \,(R(G_{pqt})) = 6$,
 \item "{ii)}" $R(G_{pqt})$ is a reducible algebraic variety when at least
 one of the absolute values of $p, q, t$ is strictly larger than two,
  \item "{iii)}" There exists an infinite set $S_2$ of groups
  associated with an infinite set $S$  of 3-tuples,  $S\subset \Bbb Z
  \times \Bbb Z \times \Bbb Z$,  having the property that if $(p,q,t)$,
  and $(p',q',t')$ are in $S$, then  $R(G_{pqt}) \ncong R(G_{p'q't'})$
 if  $p \ne p'$,\, and  $q \ne q'$,\, and $t\ne t' .$
 \endroster
 \endproclaim
 \medskip
 \demo {Proof i)}
 Using the notation of Theorem .1a let
 $P_+ = \{\rho|\rho \in R(F_2), and\, \rho(x^p_1x^q_2)=I \}$,
 $P_- = \{\rho |\rho \in R(F_2), and \, \rho (x^p_1 x^q_2)=-I \}$.
 Now by the Theorem .1a
 $$ Dim\,(R(G_{pqt}))=
 Max\, \{ Dim\,(P_+) + 2 f(t), Dim\,(P_-)+2, 3n\}, \tag 3-1.2$$
 where $f$ is the integer  valued  functions  in  the statement of
 Theorem .1. Using Theorem .1 and Lemma .1 it can be readily
 shown that
 $Dim(P_+)=4$.
 Now it is necessary to compute
 $Dim\,(P_-)$. This of course follows from
 the above lemma. Thus  $Dim(P_-)=3 \, \text {or}\, 4$,
 depending on the absolute values of $p$ and $q$. Thus, by Theorem .1a-b
 it follows that
 $$\align Dim (R(G_{pqt})) &=
 Max \{Dim(P_+) +2 f(t), Dim(P_-) +2, 6 \} \le  7  \\
 &=Max \,\{ 4+ 2f(t), (3+2) \, \text {or} \, (4+2), 6 \}=
 6. \tag 3-1.3 \endalign$$
 \enddemo
 \medskip
 \demo {Proof ii)}
 Without loss of generality it can be assumed that $t$ is
 the  integer whose
 absolute value is greater than $2$. Then in (3-1.3 ) it
 follows that since
 $t$ is not two,
 $Dim(P_+) +2 f(t)=6$. By Corollary .1 one obtains that
 $R(G_{pqt})$ is reducible.
 \enddemo
 \demo {Proof iii)}This proof will make use of Theorem 1.0.
 Beginning with the prime 3 list the infinite progression of
 primes thus:  3, 5, 7, 11, ...
 Using the resulting list associate an infinite
 set $\bold S'\subset \Bbb Z \times \Bbb Z \times \Bbb Z$
 consisting of 3-tuples as shown:\,\,
 $\bold S'=\{(3,5,7), (11,13,17),...\}.  $
 Each one of the three tuples in $\bold S'$  can
 be used to identify
 a group $G_{pqt}$; example \footnote{For clarity
 commas are inserted between $pqt$ when actual integers are used.}
 , to $(3,5,7)$ assign
 to the group
 $<x_1,x_2,x_3; x_1^{-3}x_2^{-5}=
 x_3^7>$.
 Notice also that there is a surjection
 from $G_{pqt}$ to the  free product of cyclics
 \,\, $\Bbb Z_p* \Bbb Z_q * \Bbb Z_t$, and
 $G_{pqt}/N_{pqt}\cong \Bbb Z_p* \Bbb Z_q * \Bbb Z_t,$
 where
 $N_{pqt}=N\{x_1^p, x_2^q,x_3^t\}$; here $N\{x_1^p, x_2^q,x_3^t\}$
 stands for the \lq \lq normal
 closure" in $G_{pqt}$ of the word set $\{x_1^p, x_2^q,x_3^t\}$.
 By imitation of the above, the set $\bold S'$ can be used to give rise
  to an infinite sequence
  $S_{i=1}^{\infty}(G_i,\frak N_i)$ of pairs $(G_i, \frak N_i)$
  consisting of a
  group $G_i$ and a corresponding normal subgroup $\frak N_i$ very much
  in the spirit
  of the sequence
  of Theorem 1.0, as follows:
  $$(G_1=G_{3,5,7},\frak N_1=N_{3,5,7}),(G_2=G_{11,13,17},
  \frak N_2=N_{11,13,17}),\cdots \tag 3-1.4 $$
 \flushpar Now using part i which has already been proven, it
 follows that
 $Dim(R(G_{pqt}))=6$ for any of the groups $G_{pqt}$
 in the sequence $S_{i=1}^{\infty}(G_i,\frak N_i)$.
 So as to employ Theorem 1.0, it remains to be shown
 that $Dim (R(G_{pqt}/N_{pqt}))=6$ also. But
 $G_{pqt}/N_{pqt}\cong \Bbb Z_p* \Bbb Z_q * \Bbb Z_t$. So all that
 needs to be shown is that $Dim (R(\Bbb Z_p* \Bbb Z_q * \Bbb
 Z_t))=6$, but this follows from Lemma .1 part ii, and  Lemma
 .6.
 Now all that remains to show is that the set
 $\tilde{S}=\{ N_6(R(G_j/\frak N_j))\,\vert \,(G_j,\frak N_j)
\in S_{i=1} ^{\infty}(G_i,\frak N_i)\}$ corresponding to (3-1.4)
contains an infinite set
 of integer points. But this will be a consequence of the easily
 deduced fact shown in \cite {L1}, or \cite {L2}: {\it (for $p\ge 3 $ odd \,\, $N_2(R(\Bbb Z_p))=\frac {p-1}
 {2}$).} Now using Lemma .1, Lemma .2, and
 Lemma .6 it
 follows that
 $$ N_6(R(\Bbb Z_p *\Bbb  Z_q * \Bbb Z_t))=
  \frac {(p-1)(q-1)(t-1)} {2^3}. \tag 3-1.5 $$
  Notice that $ N_6(R( \Bbb Z_p * \Bbb Z_q *\Bbb  Z_t))$ grows
  arbitrarily large as either
 $p$, or $q$, or $t$ grow large. The conditions for Theorem 1.0
 are met by (3-1.4). So one is guaranteed a
 subsequence $S^2$ of
 $S_{i=1}^{\infty}(G_i,\frak N_i)$
 with the property that given two different pairs
$(G_j,\frak N_j),\, (G_k,\frak N_k)$ in $S^2$ then
$N_6(R(G_j))\neq N_6(R(G_k))$, and a sequence $S_2$ of groups
meeting the criteria stipulated by part iii. The existence of
$S_2$ trivially leads to the existence of $S$. The proof is
complete.

 \enddemo
 \medskip
 Now some additional  remarks concerning  the  groups  in  (1-1)  are
 needed.
 \medskip
 \remark {Observation 2}
 As mentioned earlier, the groups $G_{p_1\dots p_n}$ in (1-1) are isomorphic to the groups
 $G_{-p_1\dots -p_n}$.  Thus one can assume each $p_i$ in (1-1) is
 negative as this has no consequence on the invariant $R(G)$.
 \endremark
  \medskip
 The solutions over
  $SL_2 \Bbb C$ to the equation in n-variables obtained from
 the relation of (1-1) by replacing $p_i$ with $-p_i$:
 $$a_1^{-p_1}a_2^{-p_2}\dots a_{n-1}^{-p_{n-1}}a_n ^{-p_n}=I \tag 3-1.6$$
 give rise to an affine variety in $\Bbb C^{4n}$.
 The equation in (3-1.6) can be thought of as the  relation  of
  a  one-relator
 group generated by $\{ a_1,\cdots ,a_n\}$. By Observation 2, the  solutions
 to  the  equation  (3-1.6)  gives  an  algebraic  variety  isomorphic
 to the representation variety of a group as in (1-1).
 In fact, the solutions to the equation (3-1.6) are contained in  the
 union of solutions over
  $SL_2 \Bbb C$ to the two matrix equations in $n$
 variables obtained
 from
$$a_1^{-p_1}a_2^{-p_2}\dots a_{n-1}^{-p_{n-1}}=\pm a_n ^{p_n}. \tag 3-1.7 $$
\remark { Notation }
 Denote by $\Lambda _n$ the union of solutions  over
  $SL_2 \Bbb C$  to  the  two
 matrix  equations in the sense of Observation 1 in Section One,
 obtained from (3-1.7).
 By $F_{n-1}$,  denote the free group
 generated by $\{a_1,a_2,...,a_{n-1}\}$.
 \endremark
 \medskip
  Recall that
 $R(F_{n-1})$  can  be  thought  of  as
 consisting of all  $(n-1)$-tuples of $2 \times 2$ matrices from
  $SL_2 \Bbb C$.
 Each one of the $(n-1)$-tuples is a point of $R(F_{n-1})$ and
 thus
 can be  denoted by $(m_1,m_2,...,m_{n-1})$. Again,
 the same holds true for $\Lambda_n$;  each  solution  to  the
 two matrix equations in (3-1.7) is an $n$-tuple,
 $(m_1,m_2,...,m_n)$ of $2 \times 2$  matrices
 from $SL_2 \Bbb C$.
  \medskip
 \remark {Notation}
  By $W$ denote the word $a_1^{-p_1}a_2^{-p_2}\dots
 a_{n-1}^{-p_{n-1}}$  in $F_{n-1}$, and also in the left side
  of the equation
 (3-1.7). By $\Phi$ refer to the projection  map
 $$\Phi :\Lambda_n \, \rightarrow \,R(F_{n-1}) \tag 3-1.8$$
 given by $ \Phi(m_1,m_2,\dots,m_n)=(m_1,m_2,\dots,m_{n-1}).$
 \endremark
 \medskip
 Observe that $\Lambda_n $ is the union of the two algebraic
  varieties (3-1.9)
 and (3-2.1):
 $$  \{ (\rho(a_1 ),\rho(a_2), \dots, \rho(a_{n-1}), \sigma ) \vert
 \rho \in R(F_{n-1}), \sigma \in \Omega (p_n, +\rho(W))\} \tag 3-1.9 $$
 $$  \{ (\rho(a_1 ),\rho(a_2), \dots,\rho(a_{n-1}), \sigma ) \vert
 \rho \in R(F_{n-1}), \sigma \in \Omega (p_n, - \rho(W))\}. \tag 3-2.1 $$
 \medskip
 \proclaim { Proposition 3}
 The map $\Phi$ in (3-1.8)  maps $\Lambda _n$  onto  $R(F_{n-1})$.
 \endproclaim
 \medskip
 \demo {Proof}
 One only needs to consider the case when $p_n$ is even, since
 if it is odd, then the  map is onto; see \cite {L1}.
 Suppose $p_n$ is
 even; then if
 $\rho(W)\in \,Orb(B)$, where
 by $Orb(B)$  is meant the orbit
 under $SL_2  \Bbb C$ conjugation of the
   matrix $B=
 \left(\smallmatrix  -1 & 1 \\ 0 & -1 \endsmallmatrix\right)$,
 one obtains  in  (3-1.9) that $\Omega(p_n, \rho(W))= \emptyset$.
 However,  for
 this very same choice of $\rho$,  $-\rho(W)$ does not lie in $Orb(B)$.
 Consequently in (3-2.1) the set
 $\Omega(p_n,-\rho(W))\ne \emptyset$. Thus the map  $\Phi$  is
 always onto
 since $\Lambda_n$ is the union of the algebraic varieties given  in
 (3-1.9)
 and (3-2.1).
 \enddemo
 \medskip
 Next a proposition is  proven whose  direct corollary is
 Theorem 1.2.
 \medskip
 \proclaim {Proposition 4}
 For $n \ge 3,\, Dim\,(\Lambda_n)  = 3(n-1)$.
 \endproclaim
 \medskip
 \demo { Proof}
 The proof will proceed  by  induction  on  $n \ge  3$.
  This  has  already
 been  established in Lemma .5 and Theorem 1.1 for $n = 3$. Assume
 the
 theorem
 true for $n-1$, where $n >3$. It will be then shown  true  for  $n$.
 Notice that
 $ \Lambda_n=\{(\rho(a_1),\rho(a_2),\dots,\rho(a_{n-1}),\sigma )|
 \rho \in R(F_{n-1}), \sigma \in \Omega(p_n,\pm \rho(W)) \}. $
 Given any representation
 $\rho= (m_1,...,m_{n-1}) \in R(F_{n-1})$  such that
 $ m_1^{-p_1} m_2^{-p_2}\dots  m_{n-1}^{-p_{n-1}}=\pm I;$
 then it is true that
 $Dim\,( \Phi ^{-1}(\rho)) = 2$, by Lemma .1.
 Thus the map
 $\Phi$ in general fails to have  finite  fibre; so in order to
 employ Proposition 2 care is needed.
 Let
  $U  =  \{(m_1,m_2,...,m_n)|(m_1,m_2,...,m_n) \in   \Lambda_n
  \,\, \text{and} \,\,
  m_1 ^{-p_1} \dots m_{n-1} ^{-p_{n-1}}= \pm I \}.$ \,\,
 $U$ is a sub-variety of $\Lambda_n$.  In  fact,
 $U$ is a proper sub-variety of
 $ \frak W = (\Lambda_{n-1}) \times \Omega (p_n,\pm I ),  $
 with  $Dim(\frak W)=Dim(U)$, where $\Lambda_{n-1}$
 is the union of algebraic varieties consisting of solutions over
 $SL_2 \Bbb C$
 to the two equations on $(n-1)$ variables given by
 $ a_1^{-p_1}a_2^{-p_2}\dots a_{n-2}^{-p_{n-2}}=
 \pm a_{n-1}^{p_{n-1}}. $
 Notice that by the induction hypothesis
 $$ Dim\,(\Lambda_{n-1})= 3(n-2).\tag 3-2.2$$
 Consequently, by Lemma .1,\,\,\,
  $  Dim\,(U) = 3(n-2)+2.$
 Notice that, $3(n-2)+2 < 3(n-1)$.
 Thus $U$ is a proper sub-variety of
 $\Lambda_n $ since $Dim\,(\Lambda _n) \ge 3(n-1)$, given that
 $\Lambda _n $
  contains the representation variety of a  one-relator  group
 with a presentation of deficiency
 $(n-1)$. Notice that $\{\Lambda_n- U \}$ is mapped
 by
 $\Phi$  onto $\{(R(F_{n-1}))- \Phi(U)\}$, and that $\Phi(U)$ is a
 sub-variety of
 $R(F_{n-1})$. In fact, $\Phi(U)=(\Lambda_{n-1})$. Consequently by
 (3-2.2)\,\,\,
 $ Dim\,(\Phi(U))= 3(n-2).$
 Thus $\{R(F_{n-1})- \Phi(U)\}$ is a quasi-affine variety, and since
 quasi-affine varieties are dense in an irreducible variety, it
 follows that
 $ Dim\,(\{R(F_{n-1})- \Phi(U)\})=3(n-1).$
 For every point
 $m \in \{R(F_{n-1})- \Phi(U)\}$,  the  fibre  $\Phi(m)^{-1}$
 has zero
 dimension. So present are then the needed conditions
 for Proposition 2. Hence,
 $Dim\,(\{R(F_{n-1})- \Phi(U)\}) = Dim\,(\{\Lambda_n- U\})=
 3(n-1).$
 Notice that
 $ \Lambda_n=\{ \{ \Lambda_n - U \} \cup U \}.$
 Therefore:
 $Dim(\Lambda_n)=
 Max\{Dim ( \{ \Lambda_n - U \}), Dim (U) \}
 =Max \{ 3(n-1), 3n-4 \}= 3(n-1).$
 This completes the proof.
 \enddemo

 \proclaim {Theorem 1.6} Given an integer $r\ge 2$, and an
 arbitrary integer $k\ge 1$,
 there exists a parafree group $G$ of rank $r$ and deviation one
 with $Dim(R(G))=3r$, and  $N_{3r}(R(G)) \ge k$. If $r=2$,
 the parafree group can be taken to be freely indecomposable.
 \endproclaim

 \demo {Proof}
 The infinite set of groups $S_2$ of Theorem 1.1 part iii, or Theorem 1.3,
 by virtue of
 Theorem 1.0 can be well ordered using the {\it mirc } counting function
 $N_6(V)$. So given $G_i$, $G_j$ in $S_2$ then
 either: $N_6(R(G_i))=N_6(R(G_j))$, or $N_6(R(G_i)) < N_6(R(G_j))$, or
 $N_6(R(G_i))> N_6(R(G_j))$. Now using Lemma .2, and Lemma .6
 this well ordering
 of $S_2$ can be extended using the {\it mirc } counting
 function $N_{3r}(V)$
 to the infinite set of rank $r$ deviation 1
 parafree groups $S_r$ in Theorem 1.5, where $r \geq 3$. Clearly, given
 arbitrary integer $k$,
 and $r\geq 3$, the well ordering on $S_r$ guarantees that
 there exists $G\in S_r$
 such that $N_{3r}(R(G)) \geq k$; in fact, as announced,
 an infinite number of non-isomorphic groups
 meeting the demand exist. Notice
 that each group $G$ in $S_r$ has $Dim(R(G))=3r$. Finally,
 That the groups in $S_2$ are
 freely
 indecomposable follows from Proposition 1.

 \enddemo

 \subhead {Notation}
 \endsubhead
 \flushpar $G$, a finitely generated group ({\it fg}), unless specified.
 \flushpar $SL_2\Bbb C$, group of 2 by 2 matrices of determinant one.
 \flushpar $R(G)$, the space of representations of $G$ in $SL_2\Bbb C$.
 \flushpar $Dim(R(G))$, dimension of the affine variety $R(G)$.
 \flushpar $ F_n$, the free group of rank  $n$.
 \flushpar $\gamma_n G $, the $n$-th term of the lower central series of $G$.
 \flushpar $G/ \gamma_n G$, the $n$-th term of the lower central sequence of $G$.
 \flushpar $\mu(G)$, the minimal number of generators of a group  $G$.
 \flushpar $rk(G) = \mu (G/ \gamma_2 G)$, rank of a parafree group $G$.
 \flushpar $\delta(G) = \mu(G)-rk(G)$, deviation of a parafree group $G$.
 \flushpar $W$, word in $F_n$, where $n$ depends on the context.
 \flushpar $\Sigma_n( W, k)$, solutions to the matrix equation $W=y^k$ in $SL_2\Bbb C$.
 \flushpar $\Sigma_n(-W, k)$, solutions to the equation $-W=y^k$ in $SL_2\Bbb C$.
 \flushpar $P_+=\{\rho \mid \rho \in  R(F_n),\text  {and}\,  \rho(W)=I \}.$
 \flushpar $P_-=\{\rho\mid \rho \in R(F_n),\text {and}\, \rho(W)=-I\}.$
 \flushpar $\Omega(p,M)=\{A  \mid  A\in  SL_2 \Bbb C,\,\,  A^p=M\},$ $M$
 is in $SL_2 \Bbb C$.
 \flushpar $N_c(V)$, number of $c$ dimensional maximal
 irreducible components ({\it mirc}) of an algebraic variety $V$.

  \enddocument